\documentclass[11pt,fullpage]{article}
\usepackage{amsmath}
\usepackage{amsfonts,amssymb}
\usepackage{amsthm}
\usepackage{graphics,graphicx}
\usepackage{hyperref}
\usepackage[right = 2.5cm, left=2.5cm, top = 2.5cm, bottom =2.5cm]{geometry}
\newcommand{\cE}{\mathcal E}

\newcommand{\Real}{\mathbb R}
\newcommand{\R}{\mathbb R}

\newcommand{\eps}{\epsilon}

\newcommand{\norm}[1]{\|#1\|}
\newcommand{\abs}[1]{\left\vert#1\right\vert}

\newcommand{\grad}{\nabla}

\newcommand{\Torus}{\mathbb{T}}
\newcommand{\Integer}{\mathbb{Z}}
\newcommand{\Integers}{\mathbb{Z}}

\newcommand{\Naturals}{\mathbb{N}}

\newcommand{\dd}{\, \mathrm d}

\newcommand{\ddv}{\, \mathrm dv}
\newcommand{\ddx}{\, \mathrm dx}

\newcommand{\dss}{\displaystyle}

\newcommand{\brak}[1]{\langle #1 \rangle} 

\renewcommand{\cE}{\mathcal E}

\newcommand{\cG}{\mathcal G}

\newtheorem{theorem}{Theorem}
\theoremstyle{lemma}

\theoremstyle{definition}

\theoremstyle{lemma}
\newtheorem{lemma}{Lemma}

\begin{document}

\title{A brief summary of nonlinear echoes and Landau damping} 

\author{Jacob Bedrossian\footnote{University of Maryland, College Park. Department of Mathematics and the Center for Scientific Computation and Mathematical Modeling. jacob@cscamm.umd.edu. This manuscript was completed for the Proceedings of the Journe\'es EDP 2017, based on a talk given at Journe\'es EDP 2017 in Roscoff, France. This manuscript was supported also by NSF CAREER grant DMS-1552826 and NSF RNMS \#1107444 (Ki-Net).}}
\date{\today}
\maketitle

\begin{abstract}
In this expository note we review some recent results on Landau damping in the nonlinear Vlasov equations, focusing specifically on the recent construction of nonlinear echo solutions by the author [arXiv:1605.06841] and the associated background. 
These solutions show that a straightforward extension of Mouhot and Villani's theorem on Landau damping to Sobolev spaces on $\Torus^n_x \times \Real^n_v $ is impossible and hence emphasize the subtle dependence on regularity of phase mixing problems. This expository note is specifically aimed at mathematicians who study the analysis of PDEs, but not necessarily those who work specifically on kinetic theory. However, for the sake of brevity, this review is certainly not comprehensive.
\end{abstract}

\setcounter{tocdepth}{1}
{\small\tableofcontents}
\section{Introduction and background}
We are considering the collisionless Vlasov equations: the unknown $f(t,x,v)$, called the \emph{distribution function} is a function of time, space (here we will take $x \in \Torus^n$ or $x \in \Real^n$), and velocity ($v \in \Real^n$) and satisfies: 
\begin{equation} \label{def:Vlasov}
 \left\{ 
 \begin{array}{l} 
\partial_t f + v\cdot \grad_x f + E\cdot \grad_v f = 0 \\
E = -\grad_x W \ast_x \left(\rho_f - \rho^0 \right) \\
\rho_f(t,x) = \int f(t,x,v) \ddv.
 \end{array}\right.
\end{equation}
Here $W$ is the fixed \emph{interaction potential} which is specific to the physical application and $\rho^0$ is a fixed constant whose physical relevance is explained further below.
In the physical applications we have in mind, $E = E(t,x)$ is the electric field.
The distribution function $f(t,x,v)$ is a non-negative function with the physical interpretation that given a set $A$ in phase-space $A \subset \Real^n_x \times \Real^n_v$ (or $A \subset \Torus^n_x \times \Real_v^n$), the quantity $\int_{A} f(t,x,v) \dd v \ddx$ gives the (normalized) number of particles satisfying $(x,v) \in A$. 
In the situations we discuss, the model \eqref{def:Vlasov} is mostly relevant in plasma physics, however, for some choices of $W$, this model also arises in stellar mechanics.
See e.g. \cite{Krall-Trivelpiece,BoydSanderson} for an introduction to plasma physics. 
In plasmas, $f$ can model the distribution of either electrons or ions in certain regimes and the presence of the other species is felt through $\rho^0$. 
The two most important examples of $W$ are Poisson's equation and shielded Poisson's equation: 
\begin{subequations}
\begin{align}
\widehat{W}(k) & = C_W \abs{k}^{-2}, \label{eq:Welec} \\ 
\widehat{W}(k) & = \frac{C_W}{\alpha + \abs{k}^{2}}, \label{eq:Wion}
\end{align}
\end{subequations}
where $\alpha > 0$. The latter can arise in certain regimes when modeling ions, see \cite{HanKwan11,HanKwanIacobelli14,HanKwanRousset15} and the references therein for mathematical details. 
Note the time-reversal symmetry: if $f(t,x,v)$ is a solution, then $f(-t,x,-v)$ solves the same equation backwards in time. Moreover, any re-arrangement invariant quantity is conserved -- such quantities are known as Casimirs.
There is also conservation of the total energy (if it is initially finite):
\begin{align*}
\cE = \frac{1}{2}\int \abs{v}^2 f(t,x,v) \dd v + \frac{1}{2}\int \rho W \ast \rho \dd x.
\end{align*}

Here we will be considering perturbations of a homogeneous plasma: 
\begin{align*}
f(t,x,v) = f^0(v) + h(t,x,v),
\end{align*}
with $h \in L^1_{x,v} \cap L^2_{x,v}$ (at least) and $\int h(t,x,v) \ddx \ddv =0$ and $\int f^0(v) \ddv = \rho^0$.
Note that for $x \in \Real^n_x$ we are considering localized perturbations of an infinitely extended homogeneous plasma -- basically the first question asked by physicists \cite{Vlasov-damping,Landau46,BoydSanderson}, but unfortunately it has received little attention from the mathematics community. 
The equations for $h$ are then given as
\begin{equation}  \label{def:Vlasov2}
\left\{
\begin{array}{l}
\partial_t h + v\cdot \grad_x h + E\cdot \grad_v f^0 + E\cdot \grad_v h = 0 \\
E = -\grad_x W \ast_x \rho \\
\rho(t,x) = \int h(t,x,v) \ddv \\
h(t=0,x,v) = h_{in}. 
\end{array}\right.
\end{equation}
The unknown $\rho$ is the \emph{density fluctuation}. 
For initial data satisfying, for example, $h_{in} \in H^{s;m}_{x,v}$ for an integer $m > n/2$ and $s > n/2+1$ (see appendix for notations),
it is standard to prove that the problem is locally well-posed in this class (i.e. there exists a unique solution at least on some time interval $h \in C([-T,T];H^{s;m}_{x,v})$ with $T>0$). 
Getting sharp answers to such questions is not relevant here, we will instead be concerned with understanding the behavior as $t \rightarrow \pm \infty$. 

We are concerned with understanding the stability problem, that is, we are studying solutions to \eqref{def:Vlasov2} with $h_{in}$ small in suitable spaces.
The phrase `Landau damping' refers to the behavior that $E(t) \rightarrow 0$ as $t \rightarrow \pm \infty$. Such irreversible looking behavior can seem surprising at first, however, we will see that it is not so surprising on closer inspection. 
First, we discuss the Landau damping for linear equations below in \S\ref{sec:LandauLin} and then we discuss the nonlinear problem in \S\ref{sec:nonlinearLD}. 
Further references and historical contexts are briefly discussed therein.



\section{Linear dynamics: Landau damping, phase mixing, and the Orr mechanism} \label{sec:LandauLin}
   
The linearization of Vlasov is given by: 
\begin{equation}
 \left\{ 
 \begin{array}{l} 
\partial_t h + v\cdot \grad_x h + E \cdot \grad_v f^0 = 0\\
E = -\grad_x W \ast_x \rho  \\
\rho(t,x) = \int h(t,x,v) \ddv.
 \end{array}\right.
 \end{equation}
Consider first the free transport equation (also called the free-streaming equation): 
\begin{align}
\partial_t h + v \cdot \grad_x h = 0. \label{eq:freestrm}
\end{align}
The solution, in physical and Fourier is given by: 
\begin{align}
h(t,x,v) & = h_{in}(x-tv,v) \\ 
\widehat{h}(t,k,\eta) & = \widehat{h_{in}}(t,k,\eta+kt). 
\end{align}
Therefore, the density fluctuation is given by 
\begin{align*}
\hat{\rho}(t,k) & = (2\pi)^n \widehat{h_{in}}(k,kt). 
\end{align*}
The simplest observation is that decay of the Fourier transform of the initial data implies decay in time for the density fluctuation (and more regularity $\Rightarrow$ more decay of density fluctuation):
\begin{align*}
\abs{\hat{\rho}(t,k)} & \lesssim \frac{1}{\brak{kt}^\sigma} \sup_{\eta} \abs{\brak{\eta}^\sigma \widehat{h_{in}}(k,\eta)} \\
\abs{\hat{\rho}(t,k)} & \lesssim e^{-\lambda\brak{kt}^\sigma} \sup_{\eta} \abs{e^{\lambda \brak{\eta}^\sigma}\widehat{h_{in}}(k,\eta)}. 
\end{align*}
By Sobolev embedding applied on the Fourier side, we can estimate the above supremum in frequency to velocity localization in $L^2$,
which is more amenable for doing energy estimates. 
Hence, \emph{regularity} of the initial data implies \emph{decay} of the density fluctuation.
For nonlinear problems, it is natural to determine what estimates exist which do not require a loss of regularity. 
We are not particularly interested in studying solutions with low regularity, but such estimates are important for designing quasilinear energy methods. 
In this guise, Landau damping is the same as velocity averaging (see e.g. \cite{GolseEtAl1985,GolseEtAl1988,JabinVega2004} and the references therein), though what we are discussing here is simpler than what is normally considered there.
It is a straightforward exercise in Fourier analysis and the Sobolev trace lemma (see e.g. \cite{Adams03}) to prove the following. 
\begin{lemma} \label{lem:VALD}
For all integers $m > (n-1)/2$, there holds 
\begin{align}
\norm{\rho}_{L^2_t \dot{H}^{1/2}_x} & \lesssim \norm{\brak{v}^{m}h_{in}}_{L^2_{x,v}} \\
\norm{\abs{k}^{1/2}\widehat{\rho}}_{L^\infty_k L^2_t} & \lesssim \sum_{\abs{\alpha} \leq m}\norm{D_\eta^\alpha \widehat{h_{in}}}_{L^\infty_k L^2_\eta},  
\end{align}
and at higher regularity,
\begin{subequations} \label{ineq:FreeLD}
\begin{align}
\norm{\brak{\grad_x,t\grad_x}^\sigma \rho}_{L^2_t \dot{H}^{1/2}_x} & \lesssim \norm{\brak{v}^{m}\brak{\grad_{x,v}}^\sigma h_{in}}_{L^2_{x,v}} \\
\norm{e^{\lambda\brak{\grad_x,t\grad_x}^\sigma} \rho}_{L^2_t \dot{H}^{1/2}_x} & \lesssim \norm{\brak{v}^{m}e^{\lambda\brak{\grad_{x,v}}^\sigma} h_{in}}_{L^2_{x,v}}. 
\end{align}
\end{subequations}
\end{lemma}
If one plots the evolution in phase space $(x,v) \in \Torus^n \times \Real^n$, the image is similar to milk being stirred into coffee. 
This, and related effects, are usually called \emph{phase mixing}.
To see physically what is going on, consider a very large number of non-interacting particles on a ring.
Initially, the particles might be clumped all in one spot but distributed with different speeds.  When released, they spread out around the ring and approximately distribute evenly for long stretches of time. The Landau damping in the free streaming equation is exactly this effect manifesting in the mean-field limit. 
The passage from finite to infinite dimensions in the mean-field limit is what allows Landau damping to be irreversible as $t \rightarrow \infty$: information is lost since trajectories are not compact.
This irreversibility is neither unsettling nor mysterious: it is not much different than what occurs in dispersive/wave decay \cite{TaoTextbook,KochTataruVisan14}. 
The main difference is that in dispersive equations, compactness is being lost by spreading the solution out to infinitely large scales, whereas in phase mixing problems, compactness is being lost by stirring the solution to infinitely small scales. 
Landau damping isn't any kind of damping at all, it is a type of dispersion.

In retrospect, it turns out that no discussion of Landau damping is quite complete without a discussion of the effect known to the fluid mechanics community as the \emph{Orr mechanism} \cite{Orr07} (in the plasma physics community this is sometimes called \emph{anti-phase mixing}; see e.g. \cite{SPHDDH16} and the references therein). 
Consider the initial condition 
\begin{align*} 
\widehat{h_{in}}(k,\eta) = e^{-\lambda \abs{\eta - \eta_0}}. 
\end{align*}
Then, of course
\begin{align*}
\hat{\rho}(t,k) & = (2\pi)^n e^{-\lambda \abs{kt - \eta_0}}. 
\end{align*}
Hence, the density fluctuation is exponentially localized at the specific \emph{critical time}: $t_c = \eta_0/k$. 
Specifically, this means that a small scale in the initial data (scale $\eta_0^{-1}$) will eventually be unmixed to $O(1)$ scales after a long time. 
All frequencies with $\eta_0/k > 0$ are eventually unmixed at $t_c = \eta_0/k$ to $O(1)$ scales and have their 15 minutes of fame in the electric field.  
That we can get transient growth of the density fluctuation is again not surprising: the problem is time-reversible. 
One gets a similar transient growth in dispersive equations simply by starting with dispersed initial data that is set up to re-focus at some future time 
(it is trivial to find smooth solutions to the Schr\"odinger equation $iu_t + \Delta u = 0$ with $\norm{u(0)}_\infty \leq \eps$ but such that $\sup_t \norm{u(t)}_\infty = 1$). 
In dispersive equations, one pays localization to get decay  (e.g. $\norm{e^{it\Delta}u(0)}_\infty \lesssim \abs{t}^{-n/2} \norm{u(0)}_1$) and in phase mixing, one pays regularity.
This example shows, among other things, that one cannot get Landau damping estimates in $L^\infty_t$ without paying regularity (though, just like Strichartz estimates in dispersive equations \cite{TaoTextbook,KochTataruVisan14} one can get time-averaged Landau damping; in this way, velocity averaging lemmas are probably the best analogue of Strichartz estimates here). 

\subsection{Linearized Vlasov equations}
The linearized equations are given by
\begin{equation} \label{def:linVP}
\left\{
\begin{array}{l}
\partial_t h + v \cdot \grad_x h + E(t,x) \cdot \grad_vf^0 = 0 \\
E(t,x)  = -\grad_x \widehat{W} \ast_x \rho \\ 
\rho = \int h(t,x,v) \ddv. 
\end{array}\right.
\end{equation}
There is a beautiful additional structure to these equations which makes them tractable to detailed analysis without as much technical difficulties as first guessed. 
The key is that one can deduce the following Volterra equation for the density: 
\begin{align}
\hat{\rho}(t,k) & = \widehat{h_{in}}(k,kt) - \frac{1}{(2\pi)^{n}}\int_0^t \hat{\rho}(\tau,k) \abs{k}^2 \widehat{W}(k) (t-\tau) \widehat{f^0}(k(t-\tau)) d\tau. \label{def:rhoVolt} 
\end{align}
This  re-writes the evolution entirely on the density as a Volterra equation -- only a slight generalization of an ODE, so in a sense, we have essentially explicitly diagonalized \eqref{def:linVP}. 
Such a structure is unfortunately not present in the related fluid mechanics problems  \cite{Zillinger2016,WeiZhangZhao15,WeiZhangZhao2017,BCZV17} which makes those linear problems much more difficult. However, the structure does survive (barely) when adding certain collision operators \cite{Tristani2016,B17}.
This Volterra equation was solved first by Landau \cite{Landau46} via the Laplace transform; later expanded in e.g. \cite{VKampen55,Penrose,Degond86}. See the Paley-Wiener theory \cite{Paley-Wiener} for more on using Laplace transforms to solve Volterra equations.
Such methods have been discussed several times previously \cite{BMM13} and will not be elaborated further here. 
For $(x,v) \in \Torus^n \times \Real^n$, Penrose \cite{Penrose} derived a nearly sharp stability condition for $W$ and $f^0$ which implies that the linearized Vlasov equations behave in essentially the same manner as the free transport equation.
That is, there exists $h_\infty$ such that for $s \geq 0$ and $m > n/2+1$ integers, 
\begin{align}
\lim_{t \rightarrow \infty}\norm{e^{t v\cdot \grad_x}h(t) - h_\infty}_{H^{s,m}_{x,v}} = 0; \label{eq:scatterVlasv}
\end{align}
where we denote $e^{t v\cdot \grad_x}f = f(x+tv,v)$ (the time-reversed solution to the free transport equation). 
This is in close analogy with the concept of scattering in dispersive equations. For example, for certain potentials $V$, one can show that
solutions to the linear Schr\"odinger equation $iu_t + \Delta u + V u = 0$ satisfy
\begin{align}
\lim_{t \rightarrow \infty} \norm{e^{-it\Delta} u(t) - u_\infty}_{H^s} = 0. 
\end{align}
Each of these statements quantifies the idea that the more complicated equations (linearized Vlasov and Schr\"odinger with potential) eventually converge to the free transport and free Schr\"odinger respectively. 
On $(x,v) \in \Real^n \times \Real^n$ the linearized Vlasov equations are more subtle.
If $W$ is given by Poisson interactions \eqref{eq:Welec} it was shown in \cite{glassey94,glassey95} that the linearized Vlasov equations do \emph{not} satisfy \eqref{eq:scatterVlasv}. However, if $\abs{\widehat{W}(k)} \lesssim \brak{k}^{-2}$ (e.g. for shielded Poisson \eqref{eq:Wion}), then \eqref{eq:scatterVlasv} does hold \cite{BMM16}.
Even in this latter case however, the decay of the electric field is much slower for $x \in \Real^n$ compared to $\Torus^n$, which is an issue for the stability problem for the infinitely extended plasma. 

\section{Landau damping and nonlinear echoes} \label{sec:nonlinearLD}
The main question is to determine on what time-scales the linearized predictions of Landau and others are valid for the nonlinear Vlasov equations.
That is, we are interested in studying the following informal problem:
\begin{center}
\textbf{For a fixed norm $\norm{\cdot}_X$ and all sufficiently small $\eps$, find a time-scale $t_\ast(\eps,X)$ such that the solution to the linearized Vlasov equations is a good approximation to the nonlinear Vlasov equations if $\norm{h_{in}}_X = \eps$ and $t \leq t_\ast$.}
\end{center}
Immediately after Landau's work, physicists began to ponder this question, though naturally they were not precise about the choice of $X$ -- a mistake that would lead to some confusion. 
Many physicists predicted that $t_\ast < \infty$, and, as is sometimes the case with mathematically subtle questions, the intuition was solid, but they were both right and wrong on this prediction (see e.g. discussions in \cite{Stix,MouhotVillani11,BMM13} and the references therein). 
The simplest formal calculation, found in most plasma physics texts \cite{BoydSanderson}, suggests that maybe $t_\ast \approx \eps^{-1}$ is the best one can hope for.

Work on the nonlinear Landau damping problem began in the mathematics community much later with the work of Caglioti and Maffei \cite{CagliotiMaffei98} (see also \cite{HwangVelazquez09} for improvements).
These authors exhibited real analytic solutions to the Vlasov equations with $x \in \Torus$ which displayed Landau damping as predicted by the linearized Vlasov equations.
The work is quite analogous to proving the \emph{injectivity} of wave operators in dispersive equations \cite{TaoTextbook}, e.g. proving that for nonlinear Schr\"odinger $iu_t + \Delta u = \sigma \abs{u}^p u$, for all $u_+ \in H^s$, there is a solution to NLS defined on $[T,\infty)$ (for $T$ maybe very large) such that $u(t) \approx e^{it\Delta}u_+$.
In \cite{LZ11b}, the authors characterized some non-Landau damping solutions which are close to equilibrium in  $H^s$ for $s < 3/2$, hence proving that for $X = H^s$ with $s < 3/2$, in fact, $t_\ast = 0$.
The breakthrough work of Mouhot and Villani \cite{MouhotVillani11} made major progress on the time-scale question, working with the Gevrey-$1/s$ norm:
\begin{align}
\norm{e^{\lambda \brak{\grad}^s}f}_{2} := \norm{f}_{\cG^{\lambda,s}}. 
\end{align}
The authors prove that for $X = \cG^{\lambda,s}$ for all $s \in (0,1]$ sufficiently close to one, one has $t_\ast = +\infty$. 
An alternative, simpler, proof was put forward in \cite{BMM13} which also extended the results to the full range of exponents predicted by Mouhot and Villani: $s \in (1/3,1]$.
\begin{theorem}[Mouhot/Villani \cite{MouhotVillani11} (for $s > 1/3$, see \cite{BMM13}) ]\label{thm:Main} 
Let $(x,v) \in \mathbb T^n \times \Real^n$, $f^0$ satisfy the Penrose linear stability condition, let $\frac{1}{3} < s \leq 1$, $\lambda > \lambda' > 0$, and let $m > n/2$ be an integer. Then there exists an $\epsilon_0$ such that if $h_{in}$ is mean zero, and\footnote{At least in 1D, it suffices to take $s=1/3$ and $\lambda = O(\eps^{1/3})$ \cite{B16} which shows that the requirement is less stringent than it looks as it is only relevant for $\abs{\xi} \gtrsim \eps^{-1}$ (not a coincidence if one considers the Orr mechanism).}
\begin{align*} 
\norm{e^{\lambda \brak{\grad}^s}h_{in}}_{H^{0,m}_{x,v}} = \eps < \eps_0,
\end{align*} 
then there exists a mean-zero $h_\infty$ satisfying\footnote{One can be precise in asserting that the linearized equations are leading order accurate.}
\begin{align}
    \norm{e^{\lambda' \brak{\grad}^s}\left(e^{t v\cdot \grad_x} h - h_\infty\right)}_{H^{0,m}_{x,v}} & \lesssim  \epsilon e^{-\frac{1}{2}(\lambda - \lambda^\prime)t^{s}}, \label{ineq:scatCon}\\ 
\abs{\hat{\rho}(t,k)} & \lesssim  \epsilon e^{-\lambda'\abs{kt}^s}. 
\end{align}  
\end{theorem}
Similar high regularity requirements arose also in related problems in fluid mechanics \cite{BM13,BMV14,BGM15I,BGM15II}. 
It is incorrect to argue that Gevrey regularity is completely non-physical, however, it is also incorrect to assert that this gives the whole physical picture, indeed, while perhaps some experiments might be modeled with Gevrey regular data, it is surely not the case that all physical settings of interest should be modeled thus.
In Sobolev spaces, the predictions of linearized Vlasov on $\Torus^n \times \Real^n$ are (perhaps not sharp): for $h_{in} \in H^{\sigma;m}_{x,v}$, there exists $h_\infty \in H^{s-2;m}_{x,v}$ such that for $\alpha < s-2$
\begin{align}
\norm{e^{t v \cdot \grad_x} h(t) - h_\infty}_{H^{\alpha;m}_{x,v}} & \lesssim \eps \brak{t}^{-s-2 + \alpha} \\
\norm{E(t)}_{L^2} & \lesssim \eps \brak{t}^{-\sigma}.
\end{align}
Using the techniques of \cite{BMM13}, it is relatively easy to prove that essentially the same estimates hold for the nonlinear Vlasov equations until $t_\ast \gtrsim \eps^{-1}$.
The next result, proved a few years later in \cite{B16}, shows that indeed, $t_\ast$ cannot be much larger than $\eps^{-1}$ (the proof shows that $t_\ast$ is at most just a little larger than $\abs{\log \eps} \eps^{-1}$). 
Let $\delta$ be small and $\gamma_0 \geq 2$,  
\begin{align*}
f^0(v) & = \frac{4 \pi \delta}{1+v^2} \\ 
\widehat{W}(k) & = \pm\abs{k}^{-\gamma_0}  \textbf{ or } \pm (1+\abs{k})^{-\gamma_0}.
\end{align*}
We will have to require $\delta \lesssim \eps^p$ for $p \in(0,1)$; we believe this to be purely technical.
The proof also strongly suggests that Gevrey-3 (e.g. $s=1/3$ in Theorem \ref{thm:Main}) is exactly the sharp regularity; we believe that our proof could probably be extended to prove this with some technical effort.
\begin{theorem}[\cite{B16}] \label{thm:Echo} 
Let $R \geq 1$ and $p\in(0,1)$ be arbitrary. For all $\sigma \gtrsim R$, for all $\epsilon$ sufficiently small, all $0 < \delta  \leq \epsilon^p$, there exists an $h_{in}$ satisfying\footnote{Nothing is pathological qualitatively: $h_{in}$ is also real analytic and is chosen such that $f^0 + h_{in}$ is strictly positive.} 
\begin{align}
\norm{h_{in}}_{H^{\sigma;1}_{x,v}} \leq \epsilon, \label{ineq:hinest}
\end{align}
but such that at some finite time $t_\ast = t_\ast(\epsilon,R) \gtrsim \abs{\log \eps}\eps^{-1}$: 
\begin{align}
\norm{e^{t_\ast v\cdot \grad_x} h(t_\ast)}_{H^{\sigma-R+z}} & \gtrsim t_\star^{z} \gg \epsilon^{-z}, \label{ineq:hcircTexplode} \\ 
\norm{E(t_\ast)}_{L^2} & \gtrsim t_\ast^{R-\sigma}. \label{ineq:Etstardefct}
\end{align} 
\end{theorem}
The story is a bit different on $\Real^n \times \Real^n$; see \cite{BMM16} for more discussion.

\subsection{Landau damping in the nonlinear Vlasov equations on $\Torus^n \times \Real^n$} \label{sec:LDNon}
The original proof of Theorem \ref{thm:Main} by Mouhot and Villani was obtained by a Nash-Moser-style Newton iteration \cite{MouhotVillani11}.
This is a natural way to approach the problem, but it was pointed out by H\"ormander \cite{Hormander1990} that many of the advantages of these methods
can be obtained via alternative approaches involving the paradifferential calculus of Bony \cite{Bony81}.
The intuition for this is that Nash-Moser/Newton iterations layer on corrections to the solution at successively higher frequencies by linearizing around lower frequencies. Paraproduct expansions, intuitively, linearize the evolution of higher frequencies around lower ones in a simpler and more dynamic way.
In \cite{BMM13} we introduce such an energy method for obtaining nonlinear Landau damping results\footnote{The introduction of paraproducts instead of Nash-Moser/Newton for phase mixing is originally from \cite{BM13}, but other than that, the methods of \cite{BMM13} and \cite{BM13} are different.}. The methods actually bear a little more similarity to the tricks used in quasilinear dispersive equations than to methods usually applied in fluid mechanics or kinetic theory.
The proof starts by writing (in dispersive equations $g$ is usually called the \emph{profile}), 
\begin{align}
g(t,z,v) := h(t,z+tv,v) = e^{t v\cdot \grad_x} h(t).  
\end{align}
Notice the crucial fact that
\begin{align}
\widehat{\rho}(t,k) = (2\pi)^n \widehat{g}(t,k,kt).
\end{align}
It follows what uniform regularity estimates on $g$ imply Landau damping of $E$, which in turn, implies convergence estimates such as \eqref{ineq:scatCon} by straightforward arguments. 
Define
\begin{align}
A(t,k,\eta) = e^{\lambda(t)\brak{\grad_{x,v}}^s} \brak{\grad_{x,v}}^{\sigma}, 
\end{align}
with the convention that
\begin{align}
\widehat{A f}(t,k,\eta) & = A(t,k,\eta) \hat{f}(t,k,\eta) \\
\widehat{A\rho}(t,k) & = A(t,k,kt) \hat{\rho}(t,k). 
\end{align}
We then use a bootstrap argument to propagate the following hierarchy of estimates: 
\begin{subequations} \label{ctrl:BootRes}
\begin{align}
\norm{A(t,\grad) g(t)}^2_{H^{1,m}_{x,v}} & \lesssim  \brak{t}^7 \epsilon^2 \label{ctrl:HiLocalizedB} \\
\norm{A(t,\grad) g(t)}^2_{H^{-\beta,m}_{x,v}} & \lesssim \epsilon^2 \label{ctrl:LowCommLocB} \\
\int_0^{t} \norm{A\rho(\tau)}_2^2 d\tau & \lesssim \epsilon^2. \label{ctrl:Mid} 
\end{align}
\end{subequations}
The main difficulty in closing the bootstrap is in obtaining \eqref{ctrl:Mid}; the other estimates are easy variations of standard Gevrey energy methods for transport equations (see e.g. \cite{LevermoreOliver97} or \cite{BMM13} for discussion). 
We need a self-consistent estimate on the coupled system of Volterra equations: 
\begin{align*}
\widehat{\rho}(t,k) & = \widehat{h_{in}}(k,kt) - \frac{1}{(2\pi)^n}\int_0^t \widehat{\rho}(t,k) \widehat{W}(k)\abs{k}^2(t-\tau) \widehat{f^0}(k(t-\tau)) d\tau \\ 
 & - \sum_{\ell \in \Integer^d} \int_0^t \widehat{\rho}(\tau,\ell) \widehat{W}(\ell) \ell \cdot k(t-\tau) g(\tau,k-\ell,kt-\ell \tau) d\tau.     
\end{align*}
The linear term is treated via the Laplace transform. 
By the bootstrap regularity of $g$, the most challenging contribution of the nonlinear integral comes from times when $kt \approx \ell \tau$. 
One is left with considering a toy problem essentially of the form: 
\begin{align}
\widehat{\rho}(t,k)  = Easy - \eps\sum_{\ell = k\pm 1} \int_0^t \widehat{\rho}(\tau,\ell) \widehat{W}(\ell) \ell \cdot k(t-\tau) e^{-\abs{kt-\ell \tau}} d\tau.  \label{eq:rhoToy}    
\end{align}
At $kt \approx (k+1)\tau$, $k(t-\tau) \approx t$, and hence, we see that the $\rho(\frac{kt}{k+1},k+1)$ has a very strong effect on $\rho(t,k)$.
This ``resonance'' (not a true resonance, but one can call it a ``pseudo-resonance'' \cite{Trefethen2005}) is well-known to plasma physicists -- it is called a \emph{plasma echo}, first isolated in the experiments of Malmberg, Wharton, Gould, and O'Neil in 1968 \cite{MalmbergWharton68}. 
We are not concerned with one echo, we are concerned with a cascade: $k+2 \mapsto k+1 \mapsto k \mapsto k-1 \mapsto ... \mapsto 1$.
Such cascades were indeed present in the experiments. 
This potentially creates growth due to the large coefficients; indeed, it is formally a bit analogous to solving the following linear system by back-substitution: 
\begin{equation} 
\begin{pmatrix}
1 & 2 & 0 & \cdots & \cdots & \cdots & 0 \\
0 & 1 & 2 & 0      & \cdots & \cdots & 0 \\
0 & 0 & 1 & 2      & 0      & \cdots & 0 \\
  &   &   & \cdots &        &        &   \\
0 & 0 & \cdots &\cdots &  0      & 1      & 2 \\
0 & 0 & \cdots &\cdots &  0      & 0      & 1 
\end{pmatrix}
\mathbf{x} =
\begin{pmatrix}
0 \\
0 \\
\cdots \\
\cdots \\
0 \\
1
\end{pmatrix}
;
\end{equation}
one can check that $\abs{\mathbf{x}} \approx 2^n$ where $n$ is the dimension of the matrix.
Mouhot and Villani used heuristics based on this general picture to estimate $e^{c\abs{kt}^{1/3}}$ nonlinear growth on top of the linear behavior, which prompted them to suggest that perhaps Gevrey-3 is sharp \cite{MouhotVillani11}. 
It is illuminating to see what the plasma echoes look like on the distribution function $g$, which solves (for the purpose of the toy model anyway):
\begin{align}
\partial_t \hat{g}(t,k,\eta) = Easy + \eps\sum_{\ell = k \pm 1} \widehat{\rho}(t,\ell) \widehat{W}(\ell) \ell \cdot (\eta - kt) e^{-\abs{\eta-t\ell}}. \label{eq:gToy}
\end{align}
The nonlinear growth of $\hat{g}$ coming from the echoes occurs when $\eta \approx t\ell$ -- exactly Orr's critical times! --  and so we see that the inverse cascade in $\rho$ manifests as a kind of low-to-high cascade on $g$: information is pumped from modes like $(-1, O(1))$ to mode $(k,\eta)$, then $(k-1,\eta)$, then $(k-2,\eta)$ and so forth, potentially amplifying higher regularity norms of $g$.
This is not the same kind of cascade observed in e.g. nonlinear Schr\"odinger \cite{CKSTT2010,GuardiaKaloshin2015}: we do not have large norm growth because of a small amount of `energy' moving to very high modes, we have large norm growth due to a large amount of `energy' moving to `moderately high' modes. 

\subsection{Nonlinear Echoes on $\Torus \times \Real$ }
Now, let us briefly summarize the proof of Theorem \ref{thm:Echo} in \cite{B16}. 
The purpose of this paper was to construct arbitrarily small solutions with arbitrarily long chains of plasma echoes, a logical extension of the original experiments \cite{MalmbergWharton68} and a crucial step towards understanding the regularity requirements of Theorem \ref{thm:Main}. 
The construction is slightly easier in the gravitational case, so we will show only this.
The construction is done via a modified scattering approach, sharing some similarities with the constructions in \cite{CKSTT2010,HaniEtAl2015} for the nonlinear Schr\"odinger equations, despite the fact that our nonlinear cascade mechanism and `dispersive' mechanisms are very different.
Analogous to \cite{CKSTT2010,HaniEtAl2015}, we work primarily on the unknown $e^{t v\cdot \grad_x} h(t)$ and begin by constructing a toy/reduced model which we expect to be an accurate approximation for the true solution from special initial data, and then show that there exist solutions to the toy model which exhibit the desired nonlinear growth behavior (in our case, a long chain of plasma echoes).
Our toy model is similar to that discussed above in \eqref{eq:rhoToy}, \eqref{eq:gToy}.  
Then, we prove a long-time stability estimate to show that the approximate solution is close to the true solution for long enough times. 
In our setting, the latter is more challenging than the former. 

We want to construct our solutions as corrections to linear mixing behavior, so it makes sense to approximately specify a candidate $h_\infty$ in \eqref{eq:scatterVlasv} rather than an initial condition. 
Define large parameters $t_\star = \eta_0 \gg \epsilon^{-1}$ and $k_0 \approx (\epsilon \eta_0)^{1/3}$.
Fix a time $t_{in} = \eps^{-q}$ for some $q \in (0,1)$ and set
\begin{subequations}
\begin{align}
f^L(z,v) & = 8\pi\epsilon\frac{\cos(z)}{1+v^2} \label{eq:fL} \\
f^H_{in}(z,v) & = \frac{\epsilon}{\brak{k_0,\eta_0}^{\sigma}} \frac{\cos(k_0x) \cos(\eta_0 v)}{1 + 4v^2},
\end{align}
\end{subequations}
and solve the nonlinear final time problem for $r(t,z,v) = h(t,z+tv,v)$, 
\begin{equation} \label{def:fHST}
\left\{
\begin{array}{l} \dss
\partial_t r + E(t,z+tv)(\partial_v - t\partial_z)(f^0 + r) = 0, \quad\quad t < t_{in}, \\
\rho(t,x) = \int_{\Real} r(t,z-tw,w) dw \\  
r(t_{in})  = f^L(z,v) + f^H_{in}(z,v),  
\end{array}
\right.
\end{equation}
to find the initial condition; see \cite{B16} for more information on how to perform the requisite (fairly standard) energy estimates\footnote{This implies that, while we have a very good idea about the low frequencies of the initial condition, we lose precise information about the high frequencies in the data.}.
For $t > t_{in}$ we decompose the solution (after free transport evolution) into three contributions: $f^L$, the approximate solution for the low frequencies, $f^H$ the approximate solution for the high frequencies, and $g$, the error: 
\begin{align*}
h(t,x,v) = f^L(x-tv,v) + f^H(t,x-tv,v) + g(t,x-tv,v). 
\end{align*}
We take the low frequency approximate solution to be the free transport evolution \eqref{eq:fL}.
The high frequency approximate solution solves basically Vlasov linearized around $f^0 + f^L$:
\begin{equation} \label{def:fHi}
\left\{
\begin{array}{l} \dss
\partial_t f^H + E^H(t,x+tv) \cdot \grad_v f^0 + E^H(t,x+tv) \cdot (\grad_v - t\grad_x) f^L = 0 \\ 
E^H = -\grad_x W \ast \rho^H, 
\end{array}
\right.
\end{equation}
with the minor adjustment that the term $E^L(t,x+tv)\cdot (\grad_v - t\grad_x) f^H$ has been dropped; however, for $t > t_{in}$ this term will essentially be negligible. 
The key point of partial linearization: the $\rho^H$ evolution can be written as an infinite system of Volterra equation: 
\begin{align}
\widehat{\rho^H}(t,k) & = \widehat{f^H_{in}}(k,kt) - \delta\int_{0}^t \widehat{\rho^H}(\tau,k) \widehat{W}(k)\abs{k}^2(t-\tau) e^{-\abs{k(t-\tau)}} d\tau \nonumber \\ 
 & \quad - \epsilon \sum_{\ell = k \pm 1} \int_{0}^t \widehat{\rho^H}(\tau,\ell) \widehat{W}(\ell) \ell \cdot k(t-\tau) e^{-\abs{kt - \ell \tau}} d\tau. \label{eq:rhoHVolt}
\end{align}
Hence, \eqref{def:fHi} is our toy model for the nonlinear dynamics. 
Lower bounds on $\rho^H$ can then be effectively obtained basically via back-substitution: first estimating $\rho^H(t,k_0)$, then estimating $\rho^H(t,k_{0}-1)$ and so forth.
It is crucial at this step to get extremely precise estimates on $f^H$ and $\rho^H$, specifically, we will need our upper and lower bounds to come very close to matching.
This requires some finesse, and is much easier in the case of gravitational interactions -- in the case of electrostatic interactions the linear term causes oscillations of the density which make it hard to get lower bounds. This is the main reason we take $\delta \lesssim \eps^p$ for $p \in (0,1)$, however, we suspect this is purely technical. 
See \cite{B16} for details on how to carry this out effectively. 

Although estimating $f^H$ and $\rho^H$ poses some challenges, the real work is estimating $g$, the error.
Denoting $f^E = f^H + f^L$ and $E^E = E^L + E^H$, the error solves: 
\begin{equation} \label{def:g}
\left\{
\begin{array}{l} \dss
\partial_t g + E_g(z+tv)\partial_v f^0 + E_g(z+tv)(\partial_v - t\partial_z) f^E + E^E(z+tv)(\partial_v - t\partial_z) g  \\ \quad\quad\quad + E_g(z+tv)(\partial_v - t\partial_z) g  = -\mathcal{E}, \quad t > t_{in} \\ 
\rho_g(t,x) = \int_{\Real} g(t,z-tw,w) dw, \\  
E_g(t,x) = -(\partial_x W \ast_{x} \rho_g)(t,x), \\
g(t_{in},z,v) = 0,
\end{array}
\right.
\end{equation}
where $\mathcal{E}$ (the `consistency error' of $f^E$) satisfies 
\begin{align}
\mathcal{E} & = E^L(z+tv)(\partial_v - t\partial_z)(f^L + f^H) + E^L(z+tv)  \partial_v f^0 +  E^H(z+tv)(\partial_v-t\partial_z) f^H.\label{def:cE}
\end{align}
The primary challenge is that it seems impossible to rule out that $g$ loses \emph{more} regularity than $f^H$ -- indeed, the same problematic terms arising from the linearization around $f^0 + f^L$ are present in both equations. This means $g$ must have \emph{higher} regularity than $f^H$.
There are two basic ideas for over-coming the obvious issues this introduces: (A) both $g$ and $f^H$ will be very small in weaker norms, that is, $f^H$ and $g$ are well-concentrated at high frequencies; (B) we need the regularity on $g$ to be not very far from $f^H$, this means we need sophisticated energy estimates capable of capturing very precisely the potential dynamics in $g$.
Energy methods of this general type were introduced to study mixing in fluid mechanics in \cite{BM13}, which subsequently led to many works on both the 2D and the 3D Navier-Stokes equations near Couette flow \cite{BMV14,BVW16,BGM15I,BGM15II,BGM15III}. See the recent review of \cite{BGM_Review17} for more information.
These methods introduce time-dependent Fourier multiplier norms which are adapted to the nonlinear `resonances' such as the echoes in 2D Euler near Couette flow as in the original work (nonlinear echoes have also been isolated experimentally in 2D Euler \cite{YuDriscollONeil,YuDriscoll02}).

Let us be a little more precise. The energy methods used to estimate $g$ are based on two Fourier multipliers:
\begin{align*} 
\norm{A(t,\grad) g}_{L^2}, \quad \quad \norm{B(t,\grad) g}_{L^2}. 
\end{align*}
The multiplier $B$ defines the low norm, a Gevrey-3 norm with a carefully tuned radius of regularity: for various parameters $\gamma$, $\nu(t)$, and $K$, 
\begin{align*} 
B(t,\grad) = \brak{\grad}^{\gamma} e^{\nu(t)(K\epsilon)^{1/3}\brak{\grad}^{1/3}}, 
\end{align*} 
whereas, the high norm is given by
\begin{subequations} \label{def:AG}
\begin{align}
A(t,\grad) & = \brak{\grad}^{\beta} e^{\mu(t)(K\epsilon)^{1/3}\brak{\grad}^{1/3}} G(t,\grad), \\ 
G(t,\grad) & = \left(\frac{e^{r(K\epsilon)^{1/3}\brak{\partial_v}^{1/3}}}{w(t,\partial_v)} + e^{r(K\epsilon)^{1/3}\brak{\partial_z}^{1/3}}\right). \label{def:G}
\end{align}
\end{subequations}
Here $w(t,\partial_v)$, which depends on $\eps$ and $K$, is related directly to precise upper bounds on the Volterra equation \eqref{eq:rhoHVolt}. 
The goal then is to adapt the bootstrap energy method in \cite{BMM13} using now both $A$ and $B$ based estimates:  
\begin{subequations}\label{def:bootg} 
\begin{align}
 \norm{\brak{v}\brak{\grad}^3Ag(t)}_{L^2} = & \lesssim \epsilon^2 \brak{t}^{5/2}, \label{ineq:Hig} \\ 
\norm{\brak{\grad_x,t\grad_x}^2 A \rho}_{L^2_t(I;L^2)} & \lesssim \epsilon^2, \label{ineq:HiRho} \\ 
\norm{\brak{v}Ag(t)}_{L^2} & \lesssim \epsilon^2, \label{ineq:Midg} \\ 
 \norm{\brak{\grad_x,t\grad_x}^2 B \rho}_{L^2_t(I;L^2)} & \lesssim 8\epsilon^{\sigma/5}, \label{ineq:LoRho} \\ 
\norm{\brak{v}Bg(t)}_{L^2} & \lesssim 8\epsilon^{\sigma/5}. \label{ineq:Log}
\end{align}
\end{subequations}
The design of $w$ will crucially allow us to obtain the estimate \eqref{ineq:HiRho} via the Volterra equations, however now, with the sharp Gevrey-3 with $\lambda = O(\eps^{1/3})$ property. 
The norms are tuned such that estimates of the following general form hold: 
\begin{align*}
\norm{A(q_1 q_2)}_{2} & \lesssim \norm{Aq_1}_2 \norm{\brak{\grad}^{-1}Bq_2}_2 + \norm{Aq_2}_2 \norm{\brak{\grad}^{-1}Bq_1}_2, 
\end{align*}
which implies that the small-ness of the $B$-norms will balance inverse powers of $\eps$ lost from e.g. $\norm{Af^H}_2$ (recall this has slightly lower regularity than $g$).  
That $E^H$ is so well-localized in time near the critical times and that $E^L$ is so small for $t \gtrsim t_{in}$ is also crucial.

\section{Suppressing echoes: low frequency interactions, dispersion, and collisions}
Despite the `negative' Sobolev space results of \cite{B16}, there have been several works exploring under what conditions one can obtain positive results in Sobolev spaces.
The first were those of \cite{MR3437866,FGVG}, which showed, for example, that for $\widehat{W}$ which are compactly supported in frequency,
one can obtain the Sobolev space equivalent of Theorem \ref{thm:Main}.
Indeed, each plasma echo requires a different mode of $W$, and hence compact support in frequency implies only finitely many plasma echoes are possible.
Another work was \cite{BMM16}, which considered $(x,v) \in \Real^3_x \times \Real^3_v$ with $\abs{\widehat{W}(k)} \lesssim \brak{k}^{-2}$ (necessary to ensure that the linearized problem satisfies \eqref{eq:scatterVlasv}).
It turns out that a subtle dispersive mechanism suppresses the plasma echoes for well-localized disturbances; see \cite{BMM16} for details.
This further emphasizes the differences between $x \in \Torus^n$ vs $x \in \Real^n$.
Finally, there is the work of the author \cite{B17}, which showed that Coulomb collisions (or more precisely, a simplified model of such) can suppress plasma echoes in Sobolev spaces, provided one relates the size of the data to the collision frequency.
Finding the optimal relation is analogous to the subcritical transition threshold problem in fluid mechanics; see \cite{BGM15I,BGM15III,B17,BGM_Review17} for more discussion.
Physicists have long argued that collisions should suppress nonlinear plasma echoes and make Landau damping possible for nonlinear models, and this work confirms that original intuition \cite{SuOberman1968,ONeil1968,Stix}.
Hence, we see that the works of \cite{MouhotVillani11,B16,B17} together paint an interesting picture of Landau damping in $(x,v) \in \Torus^n \times \Real^n$, confirming to some degree the intuition of physicists while at the same time showing that the nonlinear problem, despite being perturbative, is very subtle mathematically.

\appendix
\section{Notation and Fourier analysis conventions}
We denote $\brak{v} = \left( 1 + \abs{v}^2 \right)^{1/2}$ and we use the multi-index notation: given
$\alpha = (\alpha_1,\dots,\alpha_d) \in \Naturals^d$ and
$v = (v_1,\dots,v_d) \in \Real^d$ then
\begin{align*} 
  v^\alpha  = v^{\alpha_1}_1v^{\alpha_2}_2 \dots v^{\alpha_d}_d,
  \quad\quad\quad D_\eta^\alpha  = (i\partial_{\eta_1})^{\alpha_1} 
  \dots (i\partial_{\eta_d})^{\alpha_d}. 
\end{align*}
We use the short-hand
\begin{align*}
\abs{k,\eta} = \sqrt{\abs{k}^2 + \abs{\eta}^2} \quad\quad \brak{k,\eta}  = \brak{\abs{k,\eta}}. 
\end{align*}
We denote Lebesgue norms for $p,q \in [1,\infty]$ and $a,b \in \R^n$ as
\begin{align*}  
\norm{f}_{L_a^p L_b^q} = \left(\int_{\Real^n} \left(\int_{\Real^n} \abs{f(a,b)}^q \, {\rm d}b \right)^{p/q} \, {\rm d}a\right)^{1/p}  = \left(\int \left(\int \abs{f(a,b)}^q \, {\rm d}b \right)^{p/q} \, {\rm d}a\right)^{1/p}, 
\end{align*}
and analogously if $a \in \Torus^n$ or $b \in \Torus^n$. 

We use the notation $f \lesssim g$ when there exists a constant
$C > 0$ independent of the parameters of interest such that
$f \leq Cg$ (we analogously define $f \gtrsim g$).  Similarly, we use
the notation $f \approx g$ when there exists $C > 0$ such that
$C^{-1}g \leq f \leq Cg$.  We sometimes use the notation
$f \lesssim_{\alpha} g$ if we want to emphasize that the implicit
constant depends on some parameter $\alpha$.

For a function $g=g(x,v)$ with $(x,v) \in \Real^n_x \times \Real^n_v$, then we write its Fourier transform $\hat{g}(k,\eta)$ where $(k,\eta) \in \R^n \times \R^n$ with
\begin{align*} 
\hat{g}(k,\eta) & := \frac{1}{(2\pi)^{n}}\int_{\R^n \times \R^n} e^{-i x\cdot k - iv \cdot \eta} g(z,v) \dd z \dd v \\ 
g(x,v) & := \frac{1}{(2\pi)^{n}}\int_{\R^n \times \R^n} e^{i x\cdot k + iv \cdot \eta} \hat{g}(k,\eta) \dd k \dd \eta. 
\end{align*} 
We use an analogous convention for Fourier transforms to functions of
$x$ or $v$ alone.
Similarly if $(x,v) \in \Torus^n_x \times \Real^n_v$, then we write its Fourier transform $\hat{g}(k,\eta)$ where $(k,\eta) \in \Integers^n \times \R^n$ with
\begin{align*}
\hat{g}(k,\eta) := \frac{1}{(2\pi)^{n}}\int_{\R^n \times \R^n}
  e^{-i x \cdot k - iv \cdot \eta} g(z,v) \dd z \dd v, \quad 
g(z,v) := \frac{1}{(2\pi)^{n}}\sum_{\Integers^n} \int_{\R^n} e^{i z \cdot k + iv \cdot \eta} \hat{g}(k,\eta) \dd \eta. 
\end{align*}
Given a function $m \in L^\infty_{loc}$, we define the Fourier multiplier $m(\grad_{x,v})$ via
\begin{align*}
\widehat{m(\grad_{x,v}) f}(k,\eta) = m( (ik,i\eta) ) \hat{f}(k,\eta),  
\end{align*}
and analogously for functions and multipliers which depend only on space or velocity.
We denote Sobolev norms
For $\sigma \in \Real$, we denote the homogeneous and inhomogeneous Sobolev norms:
\begin{align*}
\norm{f}_{H^s(\Real^n)} & = \left( \int_{\Real^n} \brak{\eta}^{2s} \abs{\widehat{f}(\eta)}^2 \dd \eta\right)^{1/2} \\
\norm{f}_{\dot{H}^s(\Real^n)} & = \left( \int_{\Real^n} \abs{\eta}^{2s} \abs{\widehat{f}(\eta)}^2 \dd \eta\right)^{1/2}, 
\end{align*}
and we frequently drop the domain $\Real^n$ or $\Real^n_x \times \Real^n_v$ if it is clear from context. 
We make analogous definitions for functions of $(x,v) \in \Torus^n_x \times \Real^n_v$, as well as functions of either $x$ or $v$ alone. 

Finally we note the following interchangeability for localization weights in velocity and derivatives in frequency.
Let $m$ be an integer and $s \geq 0$. Then there holds 
\begin{align*}
\norm{\brak{v}^m \brak{\grad_{x,v}}^{s} f}_{L^2} & \approx_{m} \sum_{\alpha \in \Naturals^n: \abs{\alpha} \leq m}\norm{ v^\alpha \brak{\grad_{x,v}}^{s} f}_{L^2} \approx \sum_{\alpha \in \Naturals^n: \abs{\alpha} \leq m}\norm{ \brak{\grad_{x,v}}^{s} \left(v^\alpha f\right)}_{L^2} \\ 
& \approx \norm{\brak{\grad_{x,v}}^s(\brak{v}^m f)}_{L^2}. 
\end{align*}
The latter inequality follows from expanding on the Fourier side and applying Leibniz's rule. 
Accordingly, we define the norm
\begin{align*}
\norm{f}_{H^{s;m}} := \norm{\brak{v}^m \brak{\grad_{x,v}}^s f}_{L^2} \approx \norm{\brak{v}^m f}_{H^s}. 
\end{align*}

\bibliographystyle{abbrv}
\bibliography{eulereqns,dispersive}

\def\cprime{$'$}
\begin{thebibliography}{10}

\bibitem{Adams03}
R.~A. Adams and J.~J.~F. Fournier.
\newblock {\em Sobolev spaces}, volume 140 of {\em Pure and Applied Mathematics
  (Amsterdam)}.
\newblock Elsevier/Academic Press, Amsterdam, second edition, 2003.

\bibitem{B16}
J.~Bedrossian.
\newblock Nonlinear echoes and {Landau} damping with insufficient regularity.
\newblock {\em arXiv:1605.06841}, 2016.

\bibitem{B17}
J.~Bedrossian.
\newblock Suppression of plasma echoes and {Landau} damping in {Sobolev} spaces
  by weak collisions in a {Vlasov-Fokker-Planck} equation.
\newblock {\em To appear in Annals of PDE. arXiv:1704.00425}, 2017.

\bibitem{BCZV17}
J.~Bedrossian, M.~Coti~Zelati, and V.~Vicol.
\newblock Vortex axisymmetrization, inviscid damping, and vorticity depletion
  in the linearized 2d euler equations.
\newblock {\em arXiv preprint arXiv:1711.03668}, 2017.

\bibitem{BGM15I}
J.~Bedrossian, P.~Germain, and N.~Masmoudi.
\newblock Dynamics near the subcritical transition of the {3D Couette flow I:
  Below} threshold.
\newblock {\em To appear in Mem. Amer. Math. Soc., arXiv:1506.03720}, 2015.

\bibitem{BGM15II}
J.~Bedrossian, P.~Germain, and N.~Masmoudi.
\newblock Dynamics near the subcritical transition of the {3D Couette flow II:
  Above} threshold.
\newblock {\em arXiv:1506.03721}, 2015.

\bibitem{BGM15III}
J.~Bedrossian, P.~Germain, and N.~Masmoudi.
\newblock On the stability threshold for the {3D Couette} flow in {Sobolev}
  regularity.
\newblock {\em Ann. of Math.}, 157(1), 2017.

\bibitem{BGM_Review17}
J.~Bedrossian, P.~Germain, and N.~Masmoudi.
\newblock Stability of the {Couette} flow at high {Reynolds} number in {2D} and
  {3D}.
\newblock {\em arXiv:1712.02855}, 2017.

\bibitem{BM13}
J.~Bedrossian and N.~Masmoudi.
\newblock Inviscid damping and the asymptotic stability of planar shear flows
  in the {2D} {Euler} equations.
\newblock {\em Publications math{\'e}matiques de l'IH{\'E}S}, 122(1):195--300,
  2015.

\bibitem{BMM16}
J.~Bedrossian, N.~Masmoudi, and C.~Mouhot.
\newblock Landau damping in finite regularity for unconfined systems with
  screened interactions.
\newblock {\em To appear in Comm. Pure Appl. Math.}, 2016.

\bibitem{BMM13}
J.~Bedrossian, N.~Masmoudi, and C.~Mouhot.
\newblock Landau damping: paraproducts and gevrey regularity.
\newblock {\em Annals of PDE}, 2(1):1--71, 2016.

\bibitem{BMV14}
J.~Bedrossian, N.~Masmoudi, and V.~Vicol.
\newblock Enhanced dissipation and inviscid damping in the inviscid limit of
  the {Navier-Stokes} equations near the {2D Couette} flow.
\newblock {\em Arch. Rat. Mech. Anal.}, 216(3):1087--1159, 2016.

\bibitem{BVW16}
J.~Bedrossian, V.~Vicol, and F.~Wang.
\newblock The {Sobolev} stability threshold for {2D} shear flows near
  {Couette}.
\newblock {\em Journal of Nonlinear Science}, pages 1--25, 2016.

\bibitem{Bony81}
J.~Bony.
\newblock Calcul symbolique et propagation des singularit\'es pour les
  \'equations aux d\'eriv\'ees partielles non lin\'aires.
\newblock {\em Ann.Sc.E.N.S.}, 14:209--246, 1981.

\bibitem{BoydSanderson}
T.~J.~M. Boyd and J.~J. Sanderson.
\newblock {\em The physics of plasmas}.
\newblock Cambridge University Press, Cambridge, 2003.

\bibitem{CagliotiMaffei98}
E.~Caglioti and C.~Maffei.
\newblock Time asymptotics for solutions of {Vlasov-Poisson} equation in a
  circle.
\newblock {\em J. Stat. Phys.}, 92(1/2), 1998.

\bibitem{CKSTT2010}
J.~Colliander, M.~Keel, G.~Staffilani, H.~Takaoka, and T.~Tao.
\newblock Transfer of energy to high frequencies in the cubic defocusing
  nonlinear schr{\"o}dinger equation.
\newblock {\em Inventiones mathematicae}, 181(1):39--113, 2010.

\bibitem{Degond86}
P.~Degond.
\newblock Spectral theory of the linearized {Vlasov-Poisson} equation.
\newblock {\em Trans. Amer. Math. Soc.}, 294(2):435--453, 1986.

\bibitem{MR3437866}
E.~Faou and F.~Rousset.
\newblock Landau damping in {S}obolev spaces for the {V}lasov-{HMF} model.
\newblock {\em Arch. Ration. Mech. Anal.}, 219(2):887--902, 2016.

\bibitem{FGVG}
B.~Fernandez, D.~G{\'e}rard-Varet, and G.~Giacomin.
\newblock Landau damping in the kuramoto model.
\newblock In {\em Annales Henri Poincar{\'e}}, volume~17, pages 1793--1823.
  Springer, 2016.

\bibitem{glassey94}
R.~Glassey and J.~Schaeffer.
\newblock Time decay for solutions to the linearized {V}lasov equation.
\newblock {\em Transport Theory Statist. Phys.}, 23(4):411--453, 1994.

\bibitem{glassey95}
R.~Glassey and J.~Schaeffer.
\newblock On time decay rates in {L}andau damping.
\newblock {\em Comm. Part. Diff. Eqns.}, 20(3-4):647--676, 1995.

\bibitem{GolseEtAl1988}
F.~Golse, P.-L. Lions, B.~Perthame, and R.~Sentis.
\newblock Regularity of the moments of the solution of a transport equation.
\newblock {\em Journal of functional analysis}, 76(1):110--125, 1988.

\bibitem{GolseEtAl1985}
F.~Golse, B.~Perthame, and R.~Sentis.
\newblock Un r{\'e}sultat de compacit{\'e} pour les {\'e}quations de transport
  et application au calcul de la limite de la valeur propre principale d’un
  op{\'e}rateur de transport.
\newblock {\em CR Acad. Sci. Paris S{\'e}r. I Math}, 301(7):341--344, 1985.

\bibitem{GuardiaKaloshin2015}
M.~Guardia and V.~Kaloshin.
\newblock Growth of sobolev norms in the cubic defocusing nonlinear
  schr{\"o}dinger equation.
\newblock {\em Journal of the European Mathematical Society}, 17(1):71--149,
  2015.

\bibitem{HanKwan11}
D.~Han-Kwan.
\newblock Quasineutral limit of the {Vlasov-Poisson} system with massless
  electrons.
\newblock {\em Comm. Part. Diff. Eqns.}, 36(8):1385--1425, 2011.

\bibitem{HanKwanIacobelli14}
D.~Han-Kwan and M.~Iacobelli.
\newblock Quasineutral limit for vlasov--poisson via wasserstein stability
  estimates in higher dimension.
\newblock {\em Journal of Differential Equations}, 263(1):1--25, 2017.

\bibitem{HanKwanRousset15}
D.~Han-Kwan and F.~Rousset.
\newblock Quasineutral limit for {Vlasov-Poisson} with {Penrose} stable data.
\newblock {\em arXiv preprint arXiv:1508.07600}, 2015.

\bibitem{HaniEtAl2015}
Z.~Hani, B.~Pausader, N.~Tzvetkov, and N.~Visciglia.
\newblock Modified scattering for the cubic schr{\"o}dinger equation on product
  spaces and applications.
\newblock In {\em Forum of Mathematics, Pi}, volume~3. Cambridge University
  Press, 2015.

\bibitem{Hormander1990}
L.~H\"ormander.
\newblock The {N}ash-{M}oser theorem and paradifferential operators.
\newblock {\em Analysis, et cetera}, pages 429--449, 1990.

\bibitem{HwangVelazquez09}
H.~J. Hwang and J.~J.~L. Vela{\'z}quez.
\newblock On the existence of exponentially decreasing solutions of the
  nonlinear {Landau} damping problem.
\newblock {\em Indiana Univ. Math. J}, pages 2623--2660, 2009.

\bibitem{JabinVega2004}
P.-E. Jabin and L.~Vega.
\newblock A real space method for averaging lemmas.
\newblock {\em Journal de math{\'e}matiques pures et appliqu{\'e}es},
  83(11):1309--1351, 2004.

\bibitem{KochTataruVisan14}
H.~Koch, D.~Tataru, and M.~Visan.
\newblock Dispersive equations and nonlinear waves.
\newblock In {\em Oberwolfach Seminars}, volume~45. Springer, 2014.

\bibitem{Krall-Trivelpiece}
N.~Krall and A.~Trivelpiece.
\newblock {\em Principles of plasma physics}.
\newblock San Francisco Press, 1986.

\bibitem{Landau46}
L.~Landau.
\newblock On the vibration of the electronic plasma.
\newblock {\em J. Phys. USSR}, 10(25), 1946.

\bibitem{LevermoreOliver97}
D.~Levermore and M.~Oliver.
\newblock Analyticity of solutions for a generalized {Euler} equation.
\newblock {\em J. Diff. Eqns.}, 133:321--339, 1997.

\bibitem{LZ11b}
Z.~Lin and C.~Zeng.
\newblock Small {BGK} waves and nonlinear {L}andau damping.
\newblock {\em Comm. Math. Phys.}, 306(2):291--331, 2011.

\bibitem{MalmbergWharton68}
J.~Malmberg, C.~Wharton, C.~Gould, and T.~O'Neil.
\newblock Plasma wave echo.
\newblock {\em Phys. Rev. Lett.}, 20(3):95--97, 1968.

\bibitem{MouhotVillani11}
C.~Mouhot and C.~Villani.
\newblock On {Landau} damping.
\newblock {\em Acta Math.}, 207:29--201, 2011.

\bibitem{ONeil1968}
T.~M. O'Neil.
\newblock Effect of {Coulomb} collisions and microturbulence on the plasma wave
  echo.
\newblock {\em The Physics of Fluids}, 11(11):2420--2425, 1968.

\bibitem{Orr07}
W.~Orr.
\newblock The stability or instability of steady motions of a perfect liquid
  and of a viscous liquid, {Part I}: a perfect liquid.
\newblock {\em Proc. Royal Irish Acad. Sec. A: Math. Phys. Sci.}, 27:9--68,
  1907.

\bibitem{Paley-Wiener}
R.~E. A.~C. Paley and N.~Wiener.
\newblock {\em Fourier transforms in the complex domain}, volume~19 of {\em
  American Mathematical Society Colloquium Publications}.
\newblock American Mathematical Society, Providence, RI, 1987.
\newblock Reprint of the 1934 original.

\bibitem{Penrose}
O.~Penrose.
\newblock Electrostatic instability of a uniform non-{Maxwellian} plasma.
\newblock {\em Phys. Fluids}, 3:258--265, 1960.

\bibitem{SPHDDH16}
A.~Schekochihin, J.~Parker, E.~Highcock, P.~Dellar, W.~Dorland, and G.~Hammett.
\newblock Phase mixing versus nonlinear advection in drift-kinetic plasma
  turbulence.
\newblock {\em Journal of Plasma Physics}, 82(2), 2016.

\bibitem{Stix}
T.~Stix.
\newblock {\em Waves in plasmas}.
\newblock Springer, 1992.

\bibitem{SuOberman1968}
C.~Su and C.~Oberman.
\newblock Collisional damping of a plasma echo.
\newblock {\em Physical Review Letters}, 20(9):427, 1968.

\bibitem{TaoTextbook}
T.~Tao.
\newblock Nonlinear dispersive equations.
\newblock {\em {\textup{CBMS Regional Conference Series in Mathematics}}}, 106,
  2006.

\bibitem{Trefethen2005}
L.~N. Trefethen and M.~Embree.
\newblock {\em Spectra and pseudospectra: the behavior of nonnormal matrices
  and operators}.
\newblock Princeton University Press, 2005.

\bibitem{Tristani2016}
I.~Tristani.
\newblock Landau damping for the linearized {Vlasov Poisson} equation in a
  weakly collisional regime.
\newblock {\em arXiv:1603.07219}, 2016.

\bibitem{VKampen55}
N.~van Kampen.
\newblock On the theory of stationary waves in plasmas.
\newblock {\em Physica}, 21:949--963, 1955.

\bibitem{Vlasov-damping}
A.~A. Vlasov.
\newblock The vibrational properties of an electron gas.
\newblock {\em Zh. Eksp. Teor. Fiz.}, 291(8), 1938.
\newblock In russian, translation in english in {\em Soviet Physics Uspekhi},
  vol. 93 Nos. 3 and 4, 1968.

\bibitem{WeiZhangZhao15}
D.~Wei, Z.~Zhang, and W.~Zhao.
\newblock Linear inviscid damping for a class of monotone shear flow in sobolev
  spaces.
\newblock {\em Communications on Pure and Applied Mathematics}, 2015.

\bibitem{WeiZhangZhao2017}
D.~Wei, Z.~Zhang, and W.~Zhao.
\newblock Linear inviscid damping and vorticity depletion for shear flows.
\newblock {\em arXiv preprint arXiv:1704.00428}, 2017.

\bibitem{YuDriscoll02}
J.~Yu and C.~Driscoll.
\newblock Diocotron wave echoes in a pure electron plasma.
\newblock {\em {IEEE} Trans. Plasma Sci.}, 30(1), 2002.

\bibitem{YuDriscollONeil}
J.~Yu, C.~Driscoll, and T.~O`Neil.
\newblock Phase mixing and echoes in a pure electron plasma.
\newblock {\em Phys. of Plasmas}, 12(055701), 2005.

\bibitem{Zillinger2016}
C.~Zillinger.
\newblock Linear inviscid damping for monotone shear flows in a finite periodic
  channel, boundary effects, blow-up and critical sobolev regularity.
\newblock {\em Archive for Rational Mechanics and Analysis}, 221(3):1449--1509,
  2016.

\end{thebibliography}

\end{document}